\documentclass[12pt]{amsart}
\usepackage{epic,eepic}
 \newlength{\picunit}               
 \newlength{\templengthhoriz}        
 \newlength{\templengthvert}         
 \newlength{\temprule}               
 \newcount{\numlines}                
\newcommand{\encapsfig}[3]{                      
    \setbox0=\hbox{
	   \setlength{\unitlength}{#2\picunit}
	    \input #1.tex }              
   
     \templengthvert=\ht0                        
     \advance \templengthvert by \dp0            
     \advance \templengthvert by  10\unitlength  
     \advance \templengthvert by 0.5ex           
     \templengthhoriz=\wd0                       
     \advance \templengthhoriz by 0.75em         
     \kern 10\unitlength                         
     \if#3l                                      
     \box0                                       
     \kern -\templengthvert                      
     \numlines=0                               
     \loop \ifdim0pt<\templengthvert \advance\numlines by 1
				  \advance\templengthvert by -\baselineskip
				  \repeat
     \hangindent=\templengthhoriz                
     \hangafter=-\numlines                       
     \else                                       
     \temprule=\textwidth                        
     \advance \temprule by -\wd0                 
     \noindent\kern\temprule\box0                

     \kern -\templengthvert                      
     \numlines=0
     \loop \ifdim0pt<\templengthvert \advance\numlines by 1
				  \advance\templengthvert by -\baselineskip
				  \repeat
     \hangindent=-\templengthhoriz               
     \hangafter=-\numlines                       
     \fi                                         
   }

\newcommand{\tpoint}[1]{\vspace{3mm}\par\refstepcounter{subsection}{\bf \thesubsection.} 
  {\em #1. ---} }
\newcommand{\epoint}[1]{\vspace{3mm}\par\refstepcounter{subsection}{\bf \thesubsection.} 
  {\em #1.} }
\newcommand{\bpoint}[1]{\vspace{3mm}\par\refstepcounter{subsection}{\bf \thesubsection.} 
  {\bf #1.} }

\newlength{\baseunit}               
\newcommand{\bpf}{\noindent {\em Proof.  }}
\newcommand{\epf}{\qed \vspace{+10pt}}

\newcommand{\Q}{\mathbb{Q}}
\newcommand{\com}{\mathbb{C}}

\newcommand{\De}{\Delta}

\newcommand{\proj}{\mathbb P}
\newcommand{\oh}{{\mathcal{O}}}
\newcommand{\cC}{{\mathcal{C}}}

\newcommand{\cI}{{\mathcal{I}}}

\newcommand{\cH}{{\mathcal{H}}}

\newcommand{\cm}{{\mathcal{M}}}

\newcommand{\al}{\alpha}

\newcommand{\Om}{\Omega}
\newcommand{\be}{\beta}
\newcommand{\ga}{\gamma}

\newcommand{\de}{\delta}
\newcommand{\si}{\sigma}

\newcommand{\fm}{\mathfrak{m}}

\newcommand{\Ext}{\operatorname{Ext}}

\newcommand{\Spec}{\operatorname{Spec}}
\newcommand{\Sing}{\operatorname{Sing}}

\newcommand{\Hbar}{\overline{\cH}}
\newcommand{\Ibar}{\overline{\cI}}
\newcommand{\mbar}{\overline{M}}
\newcommand{\cmbar}{\overline{\cm}}
\newcommand{\bom}{ \rho^* \Omega_{\proj^1} \rightarrow \Omega_C}
 
\begin{document}
\pagestyle{plain}
\title{Recursions, formulas, and graph-theoretic interpretations of ramified coverings of the sphere by surfaces of genus 0 and 1}
\author{Ravi Vakil}
\date{November 29, 1998.}
\begin{abstract}
We derive a closed-form expression for all genus 1 Hurwitz numbers,
and give a simple new graph-theoretic interpretation of Hurwitz
numbers in genus 0 and 1.  (Hurwitz numbers essentially count
irreducible genus $g$ covers of the sphere, with arbitrary specified
branching over one point, simple branching over other specified
points, and no other branching.  The problem is equivalent to counting
transitive factorisations of permutations into transpositions.)  These
results prove a conjecture of Goulden and Jackson, and extend results
of Hurwitz and many others.
\end{abstract}
\maketitle

\section{Introduction}

The problem of enumerating factorisations of a permutation $\al \in
S_d$ into transpositions is one of long-standing interest in
combinatorics, functional analysis, knot theory, geometry, and
physics.  It is essentially equivalent to counting covers of the
Riemann sphere $\proj^1_{\com}$ with branching over $\infty$ given by
$\al$, fixed simple branching at other specified points of the sphere,
and no other branching, and it suffices to count irreducible covers of
genus $g$, for all $g$.  Hurwitz gave a simple formula when $g=0$ (for
any $\al$, \cite{h}); his result was largely forgotten until recently.
(Strehl has extended Hurwitz's idea to a complete proof, \cite{s}.)

Let $l(\al)$ be the number of cycles in $\al$.  D\'{e}nes gave a
formula for the case $g=0$ and $l(\al)=1$ (\cite{d}), and Arnol'd
extended this to $g=0$, $l(\al)=2$ (\cite{a}).  The physicists
Crescimanno and Taylor solved the case when $g=0$ and $\al$ is the
identity (\cite{ct}).  Goulden and Jackson dealt with the genus 0 case
in its entirety, independently recovering Hurwitz's result
(\cite{gj}).  Other proofs have since been given (e.g. \cite{gl}).

In positive genus, B. Shapiro, M. Shapiro, and Vainshtein have given a
striking formula (\cite{ssv}) when $l(\al) = 1$ (and $g$ is anything), involving
the co-efficients of the generating function
$$
\left( \frac { \sinh x/2} {x/2} \right)^{d-1}.
$$ 
They also give formulas for $g=1$ and $l(\al)=2$.  Graber and
Pandharipande have proved recursions for $g=0$ and 1 when $\al$ is the
identity ($g=0$ due to Pandharipande, $g=1$ to Graber and
Pandharipande, \cite{gp}) using divisor theory on the moduli space of
stable maps.  Goulden, Jackson, and Vainshtein (\cite{gjv}) have
derived formulas when $g+l(\al) \leq 6$, when $g=1$, $l(\al)=6$, and
when $g=1$ and $\al$ is the identity (the latter using the recursion
of Graber and Pandharipande).  They also conjectured a general formula
when $g=1$.  Recently, Ekedahl, Lando, M. Shapiro, and Vainshtein
announced (\cite{s}) that they have computed formulas for {\em all}
Hurwitz numbers as intersections of natural classes on $\cmbar_{g,n}$, the
moduli space of $n$-pointed genus $g$ curves (see \ref{highergenus}).

In this article, we use the space of stable maps to give a closed form
expression for all genus 0 and 1 Hurwitz numbers, generalizing the
genus 1 results described above, and proving the conjecture of
Goulden, Jackson, and Vainshtein (Corollary \ref{countG}).  En route,
we interpret the genus 0 and 1 numbers as counting graphs with simple
properties (Theorem \ref{geotree}).  This idea appears to be new (even
in genus 0) and suggests promising avenues for exploration in higher
genus.

\bpoint{Outline}
We use the theory of stable maps to $\proj^1$.  If $g=0$ or 1, on the
component of the moduli stack generically parametrizing degree $d$
covers by smooth curves, the divisor corresponding to maps ramified
above a certain fixed point is linearly equivalent to a divisor
supported on the locus of maps from singular curves (the
``boundary'').  By restricting this equivalence to appropriate
one-parameter families, we obtain recursions satisfied by Hurwitz
numbers (Theorem \ref{georec}), and the recursions determine the Hurwitz numbers (given the
fact that there is one degree 1 cover of $\proj^1$).  These recursions (and initial condition) are 
also satisfied by the solution to a certain graph-counting problem.  Finally, it is straightforward to get
a closed-form solution to the graph-counting problem.  

In Section \ref{geometry}, we derive the recursions, using results
of \cite{v}.  Readers unfamiliar with the language of algebraic
geometry may prefer to skip the section, reading only Theorem
\ref{georec}.  In Section \ref{combinatorics}, we relate the Hurwitz
numbers to the graph-counting problem, and derived closed-form
formulas. In Section \ref{discussionandspeculation}, we translate the
recursions into differential equations, and speculate on connections
to others' work and to higher genus.

\bpoint{Conventions}  
If an edge of a graph has both endpoints attached to the same vertex,
we say it is a {\em loop}.  If a connected graph has $V$ vertices and
$E$ edges, call $1-V+E$ the {\em genus} of the graph.  (Thus trees are
genus 0 graphs, and connected graphs with a single cycle have genus
1.)  When we count objects (e.g. covers of $\proj^1$, or graphs with
marked edges), if the automorphism group of the object is $G$, then
the object is counted with multiplicity $\frac 1 {|G|}$.  For example,
the number of connected genus 1 graphs on two labelled vertices with no
loops is $\frac 1 2$.

A {\em labelled partition} of $d$ is a partition in which the terms are
considered distinguished.  For example, there are $\binom 7 3$ ways of
splitting the labelled partition $\al = [1^7]$ into two labelled
partitions $\be = [1^3]$ and $\ga = [1^4]$.  We use set notation for
labelled partitions (e.g. in this example, $\al = \be \coprod \ga$, $\ga =
\al \setminus \be$).  If $\al$ is a labelled partition of $d$, let
$l(a)$ be the number of terms in $\al$, and let $\al_1$, $\al_2$,
\dots, $\al_{l(\al)}$ be the terms in the partition (so
$d=\al_1+\dots+\al_{l(\al)}$).  A set of transpositions in $S_d$ is
{\em transitive} if it generates $S_d$.  If $g$ is an integer, set
$r^g_\al := d + l(\al)+2g-2$.  Let $c^g_{\al}$ be the number of
factorizations of a fixed permutation $\si \in C_d$, with cycle structure
given by $\al$, into a transitive
product of $r^g_{\al}$ transpositions.

Let $G^g_{\al}$ be the number of smooth degree $d$ covers of $\proj^1_{\com}$
with ramification above $\infty$ given by $\al$, simple branching at
$r^g_{\al}$ other fixed points, and no other branching, where the
ramification points above $\infty$ are labelled.  Then by the
Riemann-Hurwitz formula, the covering curve has genus $g$.  By a
simple argument, $G^g_{\al} = c^g_{\al} / \prod_i \al_i$.  We call the numbers
$G^g_{\al}$ {\em Hurwitz numbers}.

Another number of previous interest is a variation of these: if
$h_{\al}$ is the size of the conjugacy class of $\al$ in $S_d$, then
$C^g_{\al} h_{\al}/d!$ is the number of smooth degree $d$ covers of
$\proj^1$ with ramification above $\infty$ given by $\al$, simple
branching at $r^g_{\al}$ other fixed points, and no other branching
(and no marking of points above $\infty$).  In \cite{gjv}, this number
is denote $\mu^g_{l(\al)}(\al)$.  These
numbers are often called Hurwitz numbers as well. 

\begin{figure}
\begin{center}

	   \setlength{\unitlength}{.1\picunit}
	    \begingroup\makeatletter\ifx\SetFigFont\undefined%
\gdef\SetFigFont#1#2#3#4#5{%
  \reset@font\fontsize{#1}{#2pt}%
  \fontfamily{#3}\fontseries{#4}\fontshape{#5}%
  \selectfont}%
\fi\endgroup%
{\renewcommand{\dashlinestretch}{30}
\begin{picture}(10216,2131)(0,-10)
\put(308,1507){\blacken\ellipse{150}{150}}
\put(308,1507){\ellipse{150}{150}}
\put(1508,1507){\blacken\ellipse{150}{150}}
\put(1508,1507){\ellipse{150}{150}}
\put(908,307){\blacken\ellipse{150}{150}}
\put(908,307){\ellipse{150}{150}}
\put(908,1507){\ellipse{1800}{600}}
\put(908,307){\ellipse{600}{600}}
\path(308,1507)(1508,1507)
\path(308,1507)(1508,1507)
\put(2408,1507){\blacken\ellipse{150}{150}}
\put(2408,1507){\ellipse{150}{150}}
\put(3608,1507){\blacken\ellipse{150}{150}}
\put(3608,1507){\ellipse{150}{150}}
\put(3008,307){\blacken\ellipse{150}{150}}
\put(3008,307){\ellipse{150}{150}}
\put(3008,1507){\ellipse{1800}{600}}
\put(3008,307){\ellipse{600}{600}}
\path(2408,1507)(3608,1507)
\path(2408,1507)(3608,1507)
\put(4508,1507){\blacken\ellipse{150}{150}}
\put(4508,1507){\ellipse{150}{150}}
\put(5708,1507){\blacken\ellipse{150}{150}}
\put(5708,1507){\ellipse{150}{150}}
\put(5108,307){\blacken\ellipse{150}{150}}
\put(5108,307){\ellipse{150}{150}}
\put(5108,1507){\ellipse{1800}{600}}
\put(5108,307){\ellipse{600}{600}}
\path(4508,1507)(5708,1507)
\path(4508,1507)(5708,1507)
\put(6608,1507){\blacken\ellipse{150}{150}}
\put(6608,1507){\ellipse{150}{150}}
\put(7808,1507){\blacken\ellipse{150}{150}}
\put(7808,1507){\ellipse{150}{150}}
\put(7208,307){\blacken\ellipse{150}{150}}
\put(7208,307){\ellipse{150}{150}}
\put(7208,1507){\ellipse{1800}{600}}
\put(7208,307){\ellipse{600}{600}}
\path(6608,1507)(7808,1507)
\path(6608,1507)(7808,1507)
\put(8708,1507){\blacken\ellipse{150}{150}}
\put(8708,1507){\ellipse{150}{150}}
\put(9908,1507){\blacken\ellipse{150}{150}}
\put(9908,1507){\ellipse{150}{150}}
\put(9308,307){\blacken\ellipse{150}{150}}
\put(9308,307){\ellipse{150}{150}}
\put(9308,1507){\ellipse{1800}{600}}
\put(9308,307){\ellipse{600}{600}}
\path(8708,1507)(9908,1507)
\path(8708,1507)(9908,1507)
\put(3008.000,1507.000){\arc{1200.000}{3.1416}{6.2832}}
\put(5108.000,1507.000){\arc{1200.000}{3.1416}{6.2832}}
\put(7508.000,1207.000){\arc{1897.367}{1.8925}{3.4633}}
\put(6308.000,607.000){\arc{1897.367}{5.0341}{6.6049}}
\put(9008.000,1207.000){\arc{1897.367}{5.9614}{7.5322}}
\put(10208.000,607.000){\arc{1897.367}{2.8198}{4.3906}}
\path(308,1507)(908,307)(1508,1507)
\path(3608,1507)(3008,307)
\path(4508,1507)(5108,307)
\end{picture}
} 
\end{center}
\caption{Counting $T^1_{\al}=4$ when $\al$ is the partition $3=1+2$ (the last two both count with multiplicity 1/2)}
\label{Tfig}
\end{figure}
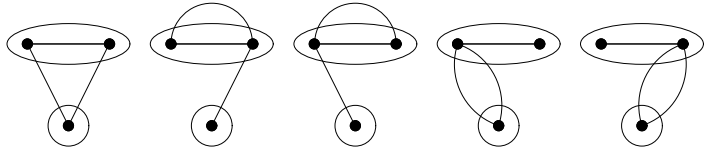

Consider $d$ labelled vertices, partitioned into subsets of size given
by a labelled partition $\al$; call these subsets {\em clumps}.   A
clump of size $i$ will be referred to as an {\em $i$-clump}.  Let
$T^g_{\al}$ be the number of connected genus $g$ graphs on these
vertices, with no loops, with a set of $\sum_i (\al_i-1)=d-l(\al)$ of
its edges that form trees on each of the clumps; we call these
$d-l(\al)$ edges the {\em edges in the clumps}.  (Thus $T^0_{\al}$
counts trees whose restriction to each of the clumps is also a tree.)
For example, if $\al$ is the partition $3=1+2$, $T^1_{\al}=4$.  This
is illustrated in Figure
\ref{Tfig}, with the clumps indicated by ovals, and the edge in
clumps indicated by drawing the edge entirely inside the corresponding 
oval.

\bpoint{Statement of results}

\tpoint{Theorem}
\label{geotree}
$$G^0_{\al} = \frac {r^0_{\al}! T^0_{\al}} {d \prod (\al_i-1)!}, \: 
G^1_{\al} = \frac {r^1_{\al}! T^1_{\al}} {12 \prod (\al_i-1)!}.$$ 

The proof is given in Section \ref{pfgeotree}.  

It is not hard to find formulas for $T^0_{\al}$, $T^1_{\al}$ (Proposition
\ref{countT}), so the above theorem gives formulas for $G^g_{\al}$
(and hence $c^g_{\al} = G^g_{\al} \prod_i \al_i$) for $g=0$, 1:

\tpoint{Corollary}
\label{countG}
{\em 
$$G^0_{\al} = \frac {r^0_{\al}!d^{l(\al)-3} \prod \al_i^{\al_i-1}
} {\prod (\al_i-1)!}, \:
G^1_{\al} = 
\frac {r^1_{\al}! d^{l(\al)-2} \prod \al_i^{\al_i-1}} {24 \prod (\al_i-1)!}
( d^2 - d - \sum_{j\geq 2} d^{2-j}(j-2)! e_j ),
$$
where $e_j$ is the $j$th symmetric polynomial in the $\al_i$. 
}

The formula for $G^0_{\al}$ is the same as that of Hurwitz, and the formula for $G^1_{\al}$ is
the conjecture of Goulden, Jackson, and Vainshtein.

\epoint{Remark}  The author has 
recently learned that Goulden and Jackson have also just proved
the above formula for $G^1_{\al}$, by purely combinatorial means
(\cite{gj98}).  The method seems unrelated.

\bpoint{Acknowledgements}
The author is grateful for discussions with A. J. de Jong, I. Goulden, M. Shapiro, and A. Vainshtein.
This project was sparked by conversations with T. Graber and R. Pandharipande.  The deformation
theory in Sections \ref{dogm}, \ref{dosm} and \ref{deftheory} was worked out
jointly with de Jong.

\section{Geometry}
\label{geometry}
Fix a labelled partition $\al$ of a positive integer $d$.  We work over the complex numbers, and
rely heavily on Sections 2--4 of \cite{v}.  All curves are assumed to
be complete.

\epoint{Background:  Stable maps to $\proj^1$}
\label{sm}
Recall that the moduli stack of stable maps $\cmbar_{g,n}(\proj^1,d)$
is a fine moduli space for degree $d$ stable maps from genus $g$ curves with $n$
labelled points to $\proj^1$.  When $g=0$, it is a smooth stack.  For
definitions and basic results, see \cite{fp}.  Let
$\cmbar_{g,n}(\proj^1,d)^+$ be the (stack-theoretic) closure in
$\cmbar_{g,n}(\proj^1,d)$ of points corresponding to maps from smooth
curves, or equivalently (from Section \ref{dosm}) the closure of
points corresponding to maps with no contracted component (i.e. where
no irreducible component of the source curve is mapped to a point).

The points of $\cmbar_{g,n}(\proj^1,d)^+$ corresponding to maps from singular curves is a union
of Weil divisors.  Such points are called {\em boundary points}.
Let $\De_0$ be the locus in $\cmbar_{g,n}(\proj^1,d)$ that is the closure
of the locus of maps of irreducible curves with one node.  If $0 \leq
i \leq g$ and $0 < j < d$, let $\De_{i,j}$ be the locus in
$\cmbar_{g,n}(\proj^1,d)$ that is the closure of maps from a reducible
curve $C_1 \cup C_2$ where $C_1$ is smooth of genus $i$ and mapping
with degree $j$, $C_2$ is smooth of genus $g-i$ and mapping with
degree $d-j$, and $C_1$ and $C_2$ meet at a node.

By \cite{v} Section 3 there is a naturally defined divisor $\be$ on
$\cmbar_g(\proj^1,d)^+$ (in the operational Chow ring) such that the
locus of maps branching above a fixed general point in $\proj^1$ lies
in class $\be [ \cmbar_g(\proj^1,d)^+]$.  If $g=0$ or 1, $\be$ is
linearly equivalent to a sum of boundary divisors (with
multiplicities).  If $g=0$, the divisor $\De_{0,j}$ appears with
multiplicity $\frac {j(d-j)} d$ (Pandharipande's relation, \cite{p}
Lemma 2.3.1, \cite{v} equation 5).  If $g=1$, the divisor $\De_0$
appears with multiplicity $\frac d {12}$, and the divisor $\De_{0,j}$
appears with multiplicity $j$ (\cite{v} Claim 4.4 and equation 6).

\epoint{Background:  Deformations of a germ of a map}
\label{dogm}

(The results of the next two sections are not surprising in the analytic category.)

We recall results about ``germs'' of maps from nodal curves to
smooth curves.  Define $\tau: \com[[z]] \rightarrow \com[[x,y]]$,
$\tau(z)=x^p+y^q$.  Let $\cC$ be the category of Artin local rings $(A,
\fm)$ over $\com$ with $A/\fm \cong \com$.  Define
the functor $F: \cC \rightarrow \text{Sets}$ as follows:

$F(A) = \{ (\de:  A[[z]] \rightarrow B, \al) \}$ up to isomorphism, where $B$ is an Artin local ring flat
over $A$,  $\al$ is an isomorphism $B \otimes_A (A / \fm) \rightarrow \com[[x,y]]/(x,y)$, and the diagram
$$
\begin{matrix} 
A[[z]] & & \stackrel {\de} \rightarrow & & B \\
\downarrow & & & & \downarrow \\
\com[[z]] & \stackrel {\tau} \rightarrow & \com[[x,y]]/(x,y) & \stackrel {\al}{\leftarrow} & B \otimes \com
\end{matrix}
$$
commutes (where the vertical arrows are restriction modulo $\fm$).  Isomorphism in this category requires
the commutativity of the obvious diagrams.

It is left to the reader to verify Schlessinger's conditions (\cite{sch}).  This functor has a hull, which
can be taken to be 
$$
R = \com [[ t,a,b_1,\dots, b_{p-1},c_1,\dots,c_{q-1}]],
$$
with $h_R \rightarrow F$ given by the ``universal curve''
\begin{equation}
\label{defspace1}
z = x^p + y^q + a + b_1 x + \dots + b_{p-1} x^{p-1} + c_1 y + \dots + c_{q-1} y^{q-1},
\end{equation}
\begin{equation}
\label{defspace2}
xy=t.
\end{equation}

Geometrically, this hull can be loosely thought of as parametrizing
deformations of the germ of a map from a node to a pointed smooth
curve (with formal co-ordinate $z$ and point $z=0$), where the node
maps to the point $z=0$, and the branches of the node ramify with
order $p$ and $q$.  The source curve is given by (\ref{defspace2}),
and the map to the pointed curve with parameter $z$ and point $z=0$ is
given by (\ref{defspace1}).  The locus where the curve remains
singular is $t=0$, which is clearly smooth (and irreducible).

Similarly, deformations of a ramification of order $p$ over a pointed
curve $z=x^p$ are given by $z=x^p + b_{p-1} x^{p-1} + \dots + b_0$.
This is well-known, and details (and the precise formulation) are left
to the reader.

\epoint{Background:  Deformations of maps to $\proj^1$}
\label{dosm}

Suppose $\rho: C \rightarrow \proj^1$ is a degree $d$ map from a nodal
curve of arithmetic genus $g$, such that no component of $C$ is
contracted.  Call formal (or analytic) neighborhoods of connected
components $A$ of $\Sing(\rho) \subset C$ {\em special loci} of
$\rho$; denote such a special locus by $(A,\rho)$.  Special loci are
(formal) neighborhoods of ramification points of $C$ or nodes of $C$.
The map $\rho$ is stable, so the functor parametrizing deformations of
the stable map $\rho$ is pro-representable by the formal neighbourhood
$X$ of the corresponding point in the moduli stack of stable maps.
The deformations are unobstructed of dimension $2d+2g-2$.  Sketch of
proof: the deformation theory of $\rho$ is controlled by $\Ext^i(\bom,
\oh_C)$.  In this case the complex $(\bom)$ is quasi-isomorphic to $(0
\rightarrow Q)$ where $Q$ is the cokernel of $\bom$, supported on the
(zero-dimensional) special loci of $\rho$.  Hence the obstruction
space $\Ext^2(\bom,\oh_C) = \Ext^2(Q,\oh_C)$ is 0, and $\Ext^2(\bom,\oh_C) = \Ext^0(Q,\oh_C)$ is 0 as well.  Then
an easy calculation using the long exact $\Ext(\cdot,\oh_C)$-sequence for
$0 \rightarrow (0 \rightarrow \Om_C) \rightarrow (\bom) \rightarrow
(\rho^* \Om_{\proj^1} \rightarrow 0) \rightarrow 0$ gives $\Ext^1(\bom,\oh) = 2d+2g-2$.
Thus $\De_0$ and $\De_{i,j}$ lie
in $\cmbar_{g,n}(\proj^1,d)^+$, and (by a quick dimension count) are
Weil divisors there.

Denote the formal schemes correspnding to the hulls of the special
loci by $X_1$, \dots, $X_n$ (whose local structure was given in
Section \ref{dogm}).  Then the natural map $X \rightarrow X_1
\times \dots \times X_n$ is an isomorphism.  A proof of this fact will appear
in \cite{ratell}; essentially it is because the sheaf $Q$ defined
in the previous paragraph is a skyscraper sheaf on the special loci,
and the restriction of $Q$ to a particular locus controls the
deformation theory of that locus.  (More generally, it will be shown
that if $\rho: C \rightarrow \proj^1$ is any stable map, perhaps with
contracted comopnents, the deformation space of $\rho$ factors into a
product of ``deformation spaces of the special loci''.)

\epoint{The stack $\Hbar^g_{\al}$}
Let $\cH^g_{\al}$ be the locus in
$\cmbar_{g}(\proj^1,d)$ corresponding to smooth covers of $\proj^1$
with ramification over $\infty$ given by $\al$, and with simple
branching over $r^g_{\al}-1$ fixed general points of $\proj^1$.  By
the Riemann-Hurwitz formula, only one ramification point is
unaccounted for.  Thus $\cH^g_{\al}$ is a one-parameter family with
one ``roaming'' simple ramification point.  Let $\Hbar^g_{\al}$ be the
(stack-theoretic) closure of $\cH^g_{\al}$ in $\cmbar_{g}(\proj^1,d)^+$.
Then $\Hbar^g_{\al}$ is a proper one-dimensional stack.

The family $\Hbar^g_{\al}$ includes points of the boundary of one of
two types.  They correspond to when the ``roaming'' ramification
hits a fixed ramification not over $\infty$, or when it hits one of
the ramifications above $\infty$.  In both cases the source curve is
either an irreducible (1-nodal) curve of geometric genus $g-1$ (i.e. where
$\Hbar^g_{\al}$ meets $\De_0$), or two smooth curves of genera adding
to $g$, joined at a node (i.e. where $\Hbar^g_{\al}$ meets some
$\De_{i,j}$).

\epoint{Multiplicity calculation}
\label{deftheory}
We can use deformation theory to compute the multiplicity
with which $\Hbar^g_{\al}$ meets the boundary divisor $\De_0$ or
$\De_{i,j}$ at a boundary point.  Let $\Hbar$ be the locus in the
hull described by (\ref{defspace1}) and (\ref{defspace2})
where the pre-image of $z=0$ remains a single point (with multiplicity
$p+q$), and the source curve is smoothed.  Suppose $I$ is the ideal of
$$
R := \com [[ t,a,b_1,\dots, b_{p-1},c_1,\dots,c_{q-1}]]
$$
defining $\Hbar$.  Multiplying (\ref{defspace1}) by $x^q$ and using $xy=t$ yields
\begin{equation}
\label{zx}
z x^q \equiv x^{p+q} + b_{p-1} x^{p+q-1} + \dots + b_1 x^{q+1} + a x^q + t^1 c_1 x^{q-1 } + \dots + 
t^{q-1} c_{q-1} x + t^q \pmod I.
\end{equation}
For convenience, let $b=\frac {b_{p-1}}{p+q}$, $c=\frac {c_{q-1}}{p+q}$, $D=\gcd(p,q)$.  When $z=0$, 
the right side of (\ref{zx}) must be a perfect $(p+q)$th power, i.e. $(x+b)^{p+q}$, so
\begin{equation} \label{check1}
b_{p-j} \equiv \binom {p+q} j b^j, \: a \equiv \binom {p+q} p b^p, \: t^j c_j \equiv \binom {p+q} {p+j} b^{p+j}, \: t^q \equiv b^{p+q} 
\pmod I,
\end{equation}
and by symmetry
\begin{equation} \label{check2}
c_{q-j} \equiv \binom {p+q} j c^j, \: a \equiv \binom {p+q} q c^q, \: t^j b_j \equiv \binom {p+q} {q+j} c^{q+j}, \: t^p \equiv c^{p+q} 
\pmod I.
\end{equation}
Note that $t c_1 = \binom {p+q} {p+1} b^{p+1}$ and $c_1 = \binom {p+q} {q-1} c^{q-1}$, so $t c^{q-1} = b^{p+1}$.

Thus  $R/I$ is generated by $b$, $c$, and $t$ with relations
$$
b^p = c^q, \: t^q = b^{p+q}, \: t^p = c^{p+q}, \: t c^{q-1} =  b^{p+1}
$$
(and possibly more).  If $D>1$, $b^p-c^q=0$ factors into $\prod_{i=1}^D
(b^{p/D} - \zeta^ic^{q/D}) = 0$, where $\zeta$ is a primitive $D$th
root of 1.  Thus $\Spec \com[b,c,t] /
(b^p-c^q,t^q-b^{p+q},t^p-c^{p+q},b^{p+1}-tc^{q-1})$ has $D$
irreducible components, with normalization parametrized by $s$, with
$b = \zeta^i s^{q/D}$, $c= s^{p/D}$, $t= \zeta^i s^{(p+q)/D}$.  (From the
last formula, each branch meets $t=0$ with multiplicity $(p+q)/D$.)
Conversely, each such branch lies in the deformations described by
(\ref{defspace1}) and (\ref{defspace2}), where the source curve is
smoothed and the pre-image of $z=0$ remains a single point, as these
branches satisfy (\ref{check1}) and (\ref{check2}).

Hence the hull of this germ of a map, keeping
ramification of order $p+q$ above $z=0$, has $D$ branches, each of
which intersect the boundary divisor $t=0$ with multiplicity
$(p+q)/D$.  Thus the intersection of $\Hbar^g_{\al}$ with the boundary at this point
is $p+q$. 

\epoint{Recursions in genus 0 and 1}

Define $\Ibar^g_{\al}$ as the closure of the locus in
$\cmbar_{g,l(\al)}(\proj^1,d)$ corresponding to smooth covers of
$\proj^1$ with ramification over $\infty$ given by $\al$, with simple
branching over $r^g_{\al}-1$ fixed general points of $\proj^1$, and
with the points over $\infty$ labelled.  In short, $\Ibar^g_{\al}$ can
be thought of as parametrizing the same maps as $\Hbar^g_{\al}$,
except the points over $\infty$ are labelled.  There is a natural
``forgetful'' morphism $\Ibar^g_{\al} \rightarrow \Hbar^g_{\al}$,
generically of degree $\prod z_{\al}(i)!$, where $z_{\al}(i)$ is the
number of times $i$ appears in $\al$.  The linear equivalences for
$\be$ in genus 0 and 1 relate the number of points on $\Hbar^g_{\al}$
where the roaming ramification maps to a fixed general point of
$\proj^1$ to the number of various boundary points (with various
multiplicities).  Each such point has $\prod z_{\al}(i)!$ pre-images
in $\Ibar^g_{\al}$.  Thus instead of counting points of
$\Hbar^g_{\al}$ in the relation, we instead count points of
$\Ibar^g_{\al}$.  (In effect, we are pulling back the relation for
$\be$ on $\Hbar^g_{\al}$ to $\Ibar^g_{\al}$.)  This will be more
convenient computationally as, for example, $\deg \be [\Ibar^g_{\al}] = G^g_{\al}$.

\tpoint{Theorem}
\label{georec}
{\em If $\al$ is a labelled partition of a positive integer $d$, then
\begin{eqnarray}
G^0_{\al} &=& (r^0_{\al}-1) \sum_{\al = \be \coprod \ga} \frac {i^2 j^2} d G^0_{\be} G^0_{\ga} 
\binom {r^0_{\al}-2} {r^0_{\be}} 
+ \sum \frac {\al_k} 2 \binom {r^0_{\al}-1} {r^0_{\be}} \frac {ij} d G^0_{\be} G^0_{\ga}, \: \: \: \rm{and}
\label{georec0} \\
G^1_{\al} &=& 2 \binom d 2 \frac d {12} ( r^1_{\al}-1) G^0_{\al} + \frac d {24} \sum \al_k G^0_{\al'} 
+ 2 (r^1_{\al}-1) \sum_{\al = \be \coprod \ga} i^2 j G^0_{\be} G^1_{\ga} \binom { r^1_{\al}-2} {r^0_{\be}}
\label{georec1}
\\
& & + \sum \al_k i G^0_{\be} G^1_{\ga} \binom { r^1_{\al}-1} {r^0_{\be}}
\nonumber
\end{eqnarray}
The first sum in (\ref{georec0}) and the second sum in (\ref{georec1}) is over
all ways of splitting $\al$ into two labelled partitions $\be$ of $i$ and $\ga$ of $j$.  The
first sum in (\ref{georec1}) is over all terms $\al_k$ of $\al$, $p+q=\al_k$, and
where $\al'$ is the labelled partition $d=\al_1 + \dots + \widehat{\al_k} + \dots +
\al_{l(\al)} + p+q$.  Similarly, the second sum in (\ref{georec0}) and the
third sum in (\ref{georec1}) is over all terms $\al_k$ of $\al$,
$p+q=\al_k$, and ways of splitting $\al'$ into two labelled partitions
$\be$ of $i$ and $\ga$ of $j$, with the $p$ in labelled partition $\be$
and the $q$ in labelled partition $\ga$.  }

Note that along with the data $G^0_{[1]} = 1$, $G^1_{[1]}=0$ (there is
one degree 1 cover of $\proj^1$, and it has genus 0), these recursions
determine $G^0_{\al}$ and $G^1_{\al}$ for all $\al$.  In the case
$\al = [1^d]$ (i.e. no ramification
over $\infty$), these are the recursions of Graber and Pandharipande
described in the introduction.

\bpf
If $g=0$, the left side of (\ref{georec0}) is $\deg
\be[\Ibar^0_{\al}]$.  By Section \ref{sm}, it can be expressed as the sum of
boundary points with certain multiplicities.  The boundary points are
of two types.

If the ``roaming'' ramification meets one of the $r^0_{\al} - 1$ fixed
ramification points over some point $P \neq \infty$, the source curve
splits into two components, one mapping with degree $i$ (say), and one
with degree $j$ --- this is a point of $\De_{0,i}$.  The ramification
over $\infty$ must be partitioned among these components, as must the
remaining $r^0_{\al} - 2$ branch points away from $\infty$.  Given
such degree $i$ and $j$ maps, there are $ij$ ways of gluing a branch
of the first curve over $P$ to a branch of the second over $P$.  There
is an additional multiplicity of $\frac {ij} d$ (from Pandharipande's
relation, Section \ref{sm}), and a multiplicity of 2 from the
multiplicity of the intersection of $\Hbar^0_{\al}$ with $\De_{0,j}$
(Section \ref{deftheory}).  Finally, we must divide by 2 because the
two components are not distinguished (any such degeneration of the
curve into $C_1 \cup C_2$ is counted twice, once
when $C_1$ corresponds to $i$, and once when it corresponds to $j$).
This gives the first term in (\ref{georec0}).

If the boundary point corresponds to when the ``roaming'' ramification
meets one of the ramifications over $\infty$ (corresponding to the
term $\al_k$, say), and the branches of the node ramify with order $p$ and $q$
($p+q=\al_k$), the source curve splits into two components, one
mapping with degree $i$, and one with degree $j$ --- this is a point of $\De_{0,i}$.  The ramification
over $\infty$ must be partitioned among these components (with the $p$
belonging to one labelled partition, and the $q$ belonging to the other), as must
the remaining $r^0_{\al} - 1$ branch points away from $\infty$.  There
is a multiplicity of $\frac {ij} d$ from Section \ref{sm}, and we must
divide by 2 because the two components are not distinguished.  By
Section \ref{deftheory}, $\Hbar^0_{\al}$ meets $\De_{0,i}$ with
multiplicity 
$\al_k$.  This gives the second term in (\ref{georec0}).

Equation (\ref{georec1}), genus 1,  is essentially the same.  

If the ``roaming'' ramification meets one of the $r^1_{\al} - 1$ fixed
ramification points over $P \neq \infty$, the source curve has a node.
First, suppose the source curve is irreducible (and genus 0) --- this
is a point of $\De_0$.  Given a map from the normalization of that
curve, there are $\binom d 2$ ways of gluing 2 different branches
above $P$ together to get a nodal curve.  There is an additional
multiplicity of $\frac d {12}$ from Section \ref{sm}, and a multiplicity
of 2 from the multiplicity of intersection of $\Hbar^1_{\al}$ with
$\De_0$ (Section \ref{deftheory}).  This gives the first term in
(\ref{georec1}).

Second, suppose the source curve splits into two components, one of
genus 0 and mapping with degree $i$ (say), and one of genus 1 and
mapping with degree $j$ --- this is a point of $\De_{0,i}$.  The ramification over $\infty$ must be
partitioned among these components, as must the remaining $r^1_{\al} -
2$ branch points away from $\infty$.  Given such degree $i$ and $j$
maps, there are $ij$ ways of gluing a branch of the first curve over
$P$ to a branch of the second over $P$.  There is an additional
multiplicity of $i$ from Section \ref{sm}, and a multiplicity of 2
from Section \ref{deftheory}.  This gives the third term in
(\ref{georec1}).

If the boundary point corresponds to when the ``roaming'' ramification
meets one of the ramifications over $\infty$ (corresponding to the
term $\al_k$, say), and the branches of the node ramify with order $p$
and $q$ ($p+q=\al_k$), suppose the source curve is irreducible ---
i.e. the boundary point lies on $\De_0$.  There is a multiplicity of
$\frac d {12}$ from Section \ref{sm}, and we must divide by 2 because
the two branches are not distinguished.  By Section \ref{deftheory},
$\Hbar^1_{\al}$ meets $\De_{0}$ with multiplicity $\al_k$.  This gives
the second term in (\ref{georec1}).

Suppose otherwise that the source curve splits into two components,
one of genus 0 mapping with degree $i$, and one of genus 1 mapping
with degree $j$ --- i.e. the boundary point lies on $\De_{0,i}$.  The
ramification over $\infty$ must be partitioned among these components
(with the $p$ belonging to one labelled partition, and the $q$
belonging to the other), as must the remaining $r^1_{\al} - 1$ branch
points away from $\infty$.  There is a multiplicity of $i$ from
Section \ref{sm}.  By Section \ref{deftheory}, $\Hbar^1_{\al}$ meets
$\De_{0,i}$ with multiplicity $\al_k$.  This gives the fourth term in
(\ref{georec1}).
\epf

\section{Combinatorics}
\label{combinatorics}

\epoint{Proof of Theorem \ref{geotree}}
\label{pfgeotree}
By inspection the result holds when $d=1$.  When $d>1$, we show that $\frac
{r^0_{\al}! T^0_{\al}} {d \prod (\al_i-1)!}$ (resp. $\frac {r^1_{\al}!
T^1_{\al}} {12 \prod (\al_i-1)!}$) satisfies the same recursion
(Theorem \ref{georec}) as $G^0_{\al}$ (resp. $G^1_{\al}$).

{\em Genus 0.}  Each tree counted by $T^0_{\al}$ has $l(\al)-1$ edges outside the clumps and $d-l(\al)$ edges inside the clumps.

The number of such trees with the choice of an edge $e$ outside the
clumps and the choice of a vertex $v$ on that edge is
$2(l(\al)-1)T^0_{\al}$.  If edge $e$ is removed, the tree breaks into
two subtrees (and each clump belongs to one of the subtrees, splitting
$\al$ into two labelled partitions $\be$ and $\ga$), with (say) $i$ and $j$
vertices respectively ($i+j=d$).  The number of ways of choosing the
subtrees, along with the vertex in each subtree to lie on $e$, is
$\sum_{\al = \be \coprod \ga} (i T^0_{\be})( j T^0_{\ga})$.  Hence
\begin{equation}
\label{pf0e1}
2(l(\al)-1) T^0_{\al} = \sum_{\al = \be \coprod \ga} i T^0_{\be} j T^0_{\ga}.
\end{equation}

The number of trees counted by $T^0_{\al}$ with the choice of an edge
$e$ {\em inside} the clump $\cC$ corresponding to $\al_k$, and a choice
of a vertex $v$ on edge $e$, is $2(\al_k-1) T^0_{\al}$.  If edge $e$
is removed, the tree breaks into subtrees, the $\al_k$ vertices in
$\cC$ are split into $p$ and $q$ ($p+q=\al_k$), and the remaining
clumps each belong to one of the subtrees as well.  The number of ways
of choosing this data is
$$
2(\al_k-1) T^0_{\al} = \sum \binom {\al_k} {p,q} (p T^0_{\be}) (q T^0_{\ga})
$$
where the sum is over all $p+q=\al_k$, and the labelled partition $\al_1 +
\dots + \widehat{\al_k} + \dots + \al_{l(\al)} +p+q$ is split into
labelled partitions $\be$ (which must contain the $p$) and $\ga$ (which must
contain the $q$).
Summing over all $\al_k$ gives
\begin{equation}
\label{pf0e2}
2(d-l(\al)) T^0_{\al} = \sum \binom {\al_k} {p,q} (p T^0_{\be}) (q T^0_{\ga})
\end{equation}
where the sum is now as described in Theorem \ref{georec}.  As
$r^0_{\al} = (d-l(\al)) + 2 (l(\al)-1)$, adding (\ref{pf0e1}) to half
(\ref{pf0e2}) gives
$$
r^0_{\al} T^0_{\al} = \sum_{\al = \be \coprod \ga} i T^0_{\be}j T^0_{\ga} + \frac 1 2 \sum \binom {\al_k} {p,q} (p T^0_{\be}) (q T^0_{\ga}).
$$
Multiplying both sides by $(r^0_{\al}-1)! / ( d \prod (\al_i-1)!)$
gives the same recursion as in Theorem \ref{georec}, with $G^0_{\de}$
replaced by $\frac {r^0_{\de}! T^0_{\de}} {d \prod (\de_i-1)!}$ for all $\de$, as desired.

{\em Genus 1.} This case is essentially the same.  Each genus 1 graph
counted by $T^1_{\al}$ has $l(\al)$ edges outside the clumps and
$d-l(\al)$ edges inside the clumps.

The number of such graphs with the choice of an edge $e$ outside the
clumps and the choice of a vertex $v$ on edge $e$ is
$2l(\al)T^1_{\al}$.  If edge $e$ is removed, then the graph either remains connected, or
breaks into two connected subgraphs.

Suppose the graph stays connected (and is hence a tree).  Then the number of ways of choosing the
tree, along with an ordered pair of distinct vertices (the endpoints of $e$), is $d(d-1)T^0_{\al}$.

Next, suppose the graph breaks into two connected subgraphs (so each
clump belongs to one of the subgraphs, splitting $\al$ into $\be$ and
$\ga$), with (say) $i$ and $j$ vertices respectively ($i+j=d$).  
One of the subgraphs is genus 0, and the other is genus 1.  The
number of ways of choosing the subgraphs, along with the vertex in
each subgraph to lie on $e$, is $ \sum_{\al = \be \coprod \ga} (i
T^0_{\be})( j T^1_{\ga}) +  (i
T^1_{\be})( j T^0_{\ga})$.  

Adding these two cases, we find:
\begin{equation}
\label{pf1e1}
2 l(\al) T^1_{\al} = d(d-1) T^0_{\al} + 2 \sum_{\al = \be \coprod \ga} (i
T^0_{\be})( j T^1_{\ga}).
\end{equation}

The number of graphs counted by $T^1_{\al}$ with the choice of an edge
$e$ {\em inside} the clump $\cC$ corresponding to $\al_k$, and a choice
of a vertex $v$ on edge $e$, is $2(\al_k-1) T^1_{\al}$.  If edge $e$
is removed, the graph either remains connected (and genus 0) or breaks into 
two subgraphs (of genus 0 and 1).   

If the graph remains connected, then edge $e$ still breaks the clump
of size $\al_k$ into two subtrees on $p$ and $q$ vertices
($p+q=\al_k$).  The number of ways of splitting the $\al_k$-clump into
a $p$-clump and a $q$-clump, then choosing the genus 0 graph, and then
choosing the vertices in the $p$-clump and $q$-clump (for endpoints
of $e$) is $\sum pq T^0_{\al'} \binom {\al_k} {p,q}$ (where $\al'$ is the
labelled partition $d=\al_1 + \dots + \widehat{\al_k} + \dots +
\al_{l(\al)} + p+q$).  

If the graph splits into two connected subgraphs, then the number of
ways of splitting the $\al_k$-clump into a $p$-clump and a $q$-clump, partitioning
the remaining clumps between $\be$ and $\ga$, choosing the endpoints of $e$ in the
$p$-clump and $q$-clump, and choosing the genus 0 and genus 1 subgraphs is
$2 \sum \binom {\al_k} {p,q} p T^0_{\be} q T^1_{\ga}$.

Adding these up over all $\al_k$ gives
\begin{equation}
\label{pf1e2}
2(d-l(\al) ) T^1_{\al} = \sum pq T^0_{\al'} \binom {\al_k} {p,q} + 
2 \sum \binom {\al_k} {p,q} p T^0_{\be} q T^1_{\ga}.
\end{equation}

As $r^1_{\al} = d+l(\al)$, adding (\ref{pf1e1}) to half
(\ref{pf1e2}) gives
$$
r^1_{\al} T^1_{\al} = d(d-1) T^0_{\al} + \frac 1 2 \sum pq T^0_{\al'} \binom {\al_k} {p,q} + 
 2 \sum_{\al = \be \coprod \ga} i T^0_{\be} j T^1_{\ga} + 
 \sum \binom {\al_k} {p,q} p T^0_{\be} q T^1_{\ga}.
$$
Multiplying both sides by $(r^1_{\al}-1)! / ( 12 \prod (\al_i-1)!)$
gives the same recursion as in Theorem \ref{georec}, with $G^g_{\de}$
replaced by $\frac {r^g_{\de}! T^g_{\de}} {12 \prod (\de_i-1)!}$ for all $\de$, and $g=0$,1 as desired.
\epf

\tpoint{Proposition} 
\label{countT}
{\em $T^0_{\al} = d^{l(\al)-2} \prod \al_i^{\al_i-1}$.  
If $e_j$ is the $j$th symmetric polynomial in the $\al_i$, 
$$
T^1_{\al} = \frac {T^0_{\al}} 2 ( d^2 - d - \sum_{j\geq 2} d^{2-j}(j-2)! e_j ).
$$
}

\bpf  The formula for $T^0_{\al}$ follows immediately from \cite{l} Ex. 4.4.

For convenience, let 
$$ S_i =
\sum_{\substack {{\al' \subset \al} \\ { l(\al')=i}}}  \left( \prod_j \al'_j \right)
\left( \sum_j \al'_j \right) = 
\sum_{\substack {{\al' \coprod \al'' =  \al} \\ { l(\al')=i}}}  \left( \prod_j \al'_j \right)
\left( d- \sum_j \al''_j \right) =  d e_i - (i+1) e_{i+1}.
$$

The number of graphs counted in $T^1_{\al}$ where the cycle passes
through only 1 clump of size $\al_k$ is $\binom {\al_k} 2 T^0_{\al}$:
one edge of the cycle is {\em outside} the clump, and there are
$T^0_{\al}$ choices for such graphs not including this edge, and
$\binom {\al_k} 2$ choices for this edges.  Summing over all $k$, we have
$\frac 1 2 T^0_{\al} (d^2-d-2e_2)$.

We next count the number of graphs counted in $T^1_{\al}$ where the
cycle passes through $i$ clumps ($i>1$).   The number of such graphs where the
cycle passes through the $i$ clumps corresponding to some $\al'
\subset \al$, $l(\al')=i$ (where $\al'' = \al \setminus \al'$) can be
counted as follows.  There are $\prod \al_j^{\al_j-2}$ ways of choosing the edges
inside the clumps (Cayley's theorem).  There are $(i-1)!/2$ ways of choosing the cyclic ordering
of the $i$ clumps in the cycle.  Then there are $\prod (\al_j')^2$ ways of choosing the
cycle itself.  Finally, by \cite{l} Ex. 4.4 there are $(\sum \al'_j)( \prod \al''_j )d^{(l(\al)-i+1) -2}$ ways
of completing the graph.  Adding these factors up over all the choices of $\al'$ gives
\begin{eqnarray*}
\prod_j \al_j^{\al_j-2} 
\sum_{\substack {{\al' \coprod \al'' =  \al} \\ { l(\al')=i}}}  \frac {(i-1)!} 2 \prod_j (\al_j')^2 \sum_j \al'_j \prod_j \al''_j d^{l(\al)-i-1}
&=& d^{l(\al)-2} \prod_j \al_j^{\al_j-1} \frac {(i-1)!} 2 d^{1-i} S_i \\
&=& T^0_{\al} \frac { (i-1)!} 2 d^{1-i} S_i
\end{eqnarray*}
Summing over all $i$, and dividing by $T^0_{\al}$,
\begin{eqnarray*}
T^1_{\al} /  T^0_{\al} &=& \frac 1 2 \left( d^2 - d - 2 e_2 + \sum_{i \geq 2} (i-1)! d^{1-i} ( d e_i - (i+1) e_{i+1} )\right) \\
&=& \frac 1 2 \left(
d^2 - d + \sum_{i \geq 2} \left( (i-1)! e_i d^{2-i}  - (i-2)! d^{2-i} i e_i \right) \right) \\
&=& \frac 1 2 \left(d^2 - d - \sum_{i \geq 2} d^{2-i} (i-2)! e_i \right)
\end{eqnarray*}
\epf

Corollary \ref{countG} follows immediately.

\section{Discussion and speculation}
\label{discussionandspeculation}
\epoint{Other recursions}
Other recursions initially seem more straightforward and natural than those of Theorem \ref{georec}.  For
example, the number of factorisations of a permutation $\si$ with given
cycle structure into transpositions $\si_1 \cdots \si_r$ can be
recursively computed by considering the possibilities for $\si_r$ (and
the possible cycle structures of $\si \si_r$), as in \cite{gj} Lemma
2.2.  However, the simplicity of the combinatorics of Section
\ref{combinatorics} suggests that the recursions of Theorem
\ref{georec} are in some way the ``right'' way to view the problem.

\epoint{Generating functions/potentials}
If $\al$ is a partition (not labelled) of $d$, define $v_{\al} =
\prod_{i \geq 1} v_i^{z_{\al}(i)}$ where $v_1$, $v_2$, \dots are formal
variables, and $z_{\al}(i)$ is the number of $i$'s in the partition $\al$. 
Recall that $h_{\al}$ was defined as the size of the conjugacy class corresponding to $\al$ in $S_d$.
Consider the following generating functions (or potential functions) for $G^0_{\al}$ and $G^1_{\al}$:
$$
F^0 = \sum_{\al \vdash d} d \frac {G^0_{\al} z^d u^{r^0_{\al}} v_{\al} h_{\al}} {d! r^0_{\al}!}, \:
F^1 = \sum_{\al \vdash d} 12 \frac {G^1_{\al} z^d u^{r^1_{\al}} v_{\al} h_{\al}} {d! r^1_{\al}!}.
$$
Both sums are over ordinary partitions (i.e. not labelled) of $d$.  ($F^0$ is similar to the generating functions $\tilde{F}$ of
\cite{gj} Lemma 2.2 and $\Phi$ of \cite{gjv} Section 3.)
  
Then Theorem
\ref{georec} can be rephrased as a differential equation, analogous to
the differential equation satisfied by the genus 0 Gromov-Witten
potential (\cite{fp}), or the differential equations
satisfied by potentials for characteristic numbers of plane curves
(\cite{ek} Section 6, \cite{v}):
\begin{eqnarray}
F^0_u &=& u z^2 ( F^0_z)^2 + \frac 1 2 \sum_{p,q} pq v_{p+q} F^0_{v_p} F^0_{v_q}, \\
F^1_u &=& u z^2 (F^0_{zz} + 2  F^0_z F^1_z) + \sum_{p,q} pq v_{p+q} \left(\frac 1 2 F^0_{v_pv_q} + F^0_{v_p}F^1_{v_q} \right).
\label{diffeq1}
\end{eqnarray}
In the genus 1 equation (\ref{diffeq1}), the second term corresponds to the third term of (\ref{georec1}) and vice versa.
As these equations do not seem especially enlightening, the details of their derivations are omitted.

\epoint{Higher genus}
\label{highergenus}

The form of Theorem \ref{geotree} is striking, and suggests immediate
generalizations to higher genus.  However, none of the obvious
extensions seem to work.  Surprising and beautiful partial results in
higher genus are already known (described in the introduction), and it
would be of interest to try to give these results a similar
graph-theoretic interpretation.

Other approaches to Hurwitz numbers also involve graph enumeration
problems.  In particular, \cite{a} and \cite{ssv} both involved
edge-ordered graphs.  It would be worthwhile to understand the
relationship between this graph-counting problem and that of this
article, as the two interpretations are advantageous in different
circumstances (here when $g \leq 1$, there when $l(\al)$ is very
small).

Another promising direction seems to be through the work of Ekedahl et
al, who express Hurwitz numbers as Hodge integrals on $\cmbar_{g,n}$.
The recursions of Theorem \ref{georec} should follow from their
analysis, as they are a consequence of the fact that the Hodge class
is linearly equivalent to a sum of boundary divisors when $g=0$ or 1.
Also, intersections on $\cmbar_{g,n}$ are naturally sums over graphs,
so the logical next step would be to try to give graph-theoretic
interpretations to higher genus Hurwitz numbers.

Graber and Pandharipande have conjectured (\cite{gp}) a recursion for genus 2 Hurwitz numbers with no ramification over $\infty$ (i.e. $\al$ is
the labelled partition $d = 1 + \dots + 1$):
\begin{eqnarray*}
G^2_{[1^d]}  &=&  d^2 \left( \frac {97}{136}  d - \frac {20}{17} \right) G^1_{[1^d]}
   + \sum_{j=1}^{d-1} G^0_{[1^j]} G^2_{[1^{d-j}]} \binom {2d} {2j-2} j(d-j)
                                   \left( -\frac{115}{17} j + 8d \right)\\
  & &  + \sum_{j=1}^{d-1} G^1_{[1^j]} G^1_{[1^{d-j}]} \binom {2d}{2j} j(d-j)
              \left( \frac{11697}{34} j(d-j) - \frac{3899}{68} d^2 \right).
\end{eqnarray*}
It is unclear why a genus 2 relation should exist (either
combinatorially or algebro-geometrically).  The relation looks as
though it is induced by a relation in the Picard group of the moduli
space, but no such relation exists.  A proof of this conjecture may
shed some light on the geometry of genus 2 pointed curves through the
work of Ekedahl et al.

The grand motivating problem behind all of these results is that of
enumerating the factorisations of a permutation $\si \in S_d$ into $r$
transpositions (not necessarily transitive) for {\em any} $d$, $\si$
and $r$, and of giving this number concrete combinatorial meaning.
One might speculate that this number would be an appropriate multiple of the number of graphs on $d$
vertices (not necessarily connected), with clumps given by $\si$, with
$d-1+g$ edges (where $g$ is given by $r=d+l(\si)+2g-2$), with some
additional structure.


\begin{thebibliography}{[SSV]}
\bibitem[A] {a} V. I. Arnol'd, {\em Topological classification of trigonometric polynomials and combinatorics of
graphs with an equal number of vertices and edges}, Functional analysis and its applications {\bf 30} no. 1 (1996) 1--17. 
\bibitem[CT] {ct}  M. Crescimanno and W. Taylor, {\em Large $N$ phases of chiral QCD${}_2$}, Nuclear Phys. B {\bf 437} no. 1 (1995) 3--24.
\bibitem[D]{d} J. D\'{e}nes, {\em The representation of a permutation as the product of a minimal number of transpositions and
its connection with the theory of graphs}, Publ. Math. Inst. Hungar. Acad. Sci. {\bf 4} (1959) 63--70.
\bibitem[EK] {ek} L. Ernstr\"{o}m and G. Kennedy, {\em Contact cohomology of the projective plane}, preprint 1998, math.AG/9703013.
\bibitem[FP] {fp} W. Fulton and R. Pandharipande, {\em Notes on stable maps and quantum cohomology}, in 
{\em Algebraic geometry Santa Cruz 1995} v. 2, J. Koll\'{a}r, R. Lazarsfeld, D. Morrison eds., A.M.S.,  Providence, 1997.
\bibitem[GL] {gl} V. V. Goryunov and S.K. Lando, {\em On enumeration of meromorphic functions of the line}, preprint.
\bibitem[GJ1] {gj} I. P. Goulden and D. M. Jackson, {\em Transitive factorisations into transpositions and holomorphic mappings
on the sphere}, Proc. of the A.M.S. {\bf 125} no. 1 (1997) 51--60.
\bibitem[GJ2] {gj98} I. P. Goulden and D. J. Jackson, manuscript
in preparation.
\bibitem[GJV]{gjv} I. P. Goulden, D. M. Jackson and A. Vainshtein, {\em The number of ramified coverings
of the sphere by the torus and surfaces of higher genera}, preprint 1998.
\bibitem[GP] {gp} T. Graber and R. Pandharipande, personal communication.
\bibitem[H] {h} A. Hurwitz, {\em \"{U}ber die Anzahl der Riemann'schen Fl\"{a}chen mit gegebenen Verzweigungspunkten}, Math. 
Ann. {\bf 55} (1902) 53--66.
\bibitem[L] {l} L. Lovasz, {\em Combinatorial Problems and Exercises} 2nd ed., Akad\'{e}miai Kiad\'{o}, Budapest, 1993. 
\bibitem[P] {p} R. Pandharipande, {\em Intersection of $\Q$-divisors on Kontsevich's moduli 
space $\mbar_{0,n}(\proj^r,d)$ and enumerative geometry}, {\em Trans. A.M.S.}, to appear.
\bibitem[Sch] {sch} M. Schlessinger, {\em Functors of Artin rings}, Trans. A.M.S. {\bf 130} (1968) 208--222.
\bibitem[S] {s} M. Shapiro, personal communication.
\bibitem[SSV] {ssv} B. Shapiro, M. Shapiro, and A. Vainshtein, {\em Ramified coverings of $S^2$ with one degenerate
branching point and enumeration of edge-ordered graphs}, Adv. in Math. Sci. {\bf 34} (1997) 219--228.
\bibitem[St] {st} V. Strehl, {\em Minimal transitive products of transpositions --- the reconstruction of a proof of
A. Hurwitz}, S\'{e}m. Lothar. Combin. {\bf 37} (1996), Art. S37c, 12 pp. (electronic, see http://cartan.u-strasbg.fr:80/\~{}slc/).
\bibitem[V1] {v} R. Vakil, {\em Recursions for characteristic numbers of genus one plane curves}, preprint 1998, submitted 
for publication; available at http//www-math.mit.edu/\~{}vakil.
\bibitem[V2] {ratell} R. Vakil, {\em The enumerative geometry of rational and elliptic curves in projective space}, manuscript in preparation;
earlier version math.AG/9709007.
\end{thebibliography}
\end{document}